\documentclass[final]{amsart}
\usepackage{pdfsync}
\usepackage{amsmath,amssymb}
\usepackage{enumerate}
\usepackage{xspace}
\usepackage{boxedminipage}
\usepackage{algorithmicx}
\usepackage[ruled]{algorithm}
\usepackage{algpseudocode}
\usepackage{listings}
\usepackage{tabularx}
\usepackage{subfigure}
\usepackage{tristan}
\usepackage[cm]{fullpage}
\renewcommand{\bi}[2]{\ensuremath{\cA\qp{#1,#2}}}

\renewcommand{\vec}[1]{\geovec{#1}}
\renewcommand{\mat}[1]{\geomat{#1}}

\renewcommand{\enorm}[2]{\ensuremath{\Norm{#1}}_{dG,{#2}}}

\newcommand{\sobg}[2]{\ensuremath{\WW^{#1,#2}_g}}
\renewcommand{\sobz}[2]{\ensuremath{\WW^{#1,#2}_0}}
\renewcommand{\sob}[2]{\ensuremath{\WW^{#1,#2}}}
\renewcommand{\leb}[1]{\ensuremath{\LL^{#1}}}

\numberwithin{equation}{section}
\setboolean{showtodo}{true}
\setboolean{shownotes}{true}
\setboolean{showchanges}{false}
\setboolean{usemathrsfs}{false}
\usepackage{tkz-graph}
\usetikzlibrary{arrows}

\author{
  Tristan Pryer
}
\address{
  Tristan Pryer
  \thanks{
    Department of Mathematics and Statistics, Whiteknights, University of Reading, Reading RG6 6AX, UK
    {\tt{T.Pryer@reading.ac.uk}}.
}}

\title[FE approximation of $\infty$-Harmonic Functions]
      {On the finite element approximation of infinity-harmonic functions}
\date{\today}

\begin{document}
\maketitle
\begin{abstract}
  In this note we show that conforming Galerkin approximations for
  $p$-harmonic functions tend to $\infty$-harmonic functions in the
  limit $p\to \infty$ and $h\to 0$, where $h$ denotes the Galerkin discretisation parameter.
\end{abstract}

\section{Introduction and the $\infty$-Laplacian}
\label{sec:introduction}

Let $\W \subset \reals^d$ be an open and bounded set. For a given
function $u : \W \to \reals$ we denote the gradient of $u$ as $\D u :
\W \to \reals^d$ and its Hessian $\Hess u:\W \to \reals^{d\times
  d}$. The $\infty$-Laplacian is the partial differential equation
(PDE)
\begin{equation}
  \label{eq:inflap}
  \Delta_\infty u 
  :=
  \frob{\qp{\D u \otimes \D u}}{\Hess u}
  =
  \sum_{i,j=1}^d \partial_i u \ \partial_j u \ \partial^2_{ij} u = 0,
\end{equation}
where ``$\otimes$'' is the tensor product between $d$-vectors and
``$:$'' the Frobenius inner product between matrices.

This problem is the prototypical example of a PDE from Calculus of
Variations in $\leb{\infty}$, arising as the analogue of the \emph{Euler--Lagrange}
equation of the functional 
\begin{equation}
  \cJ[u;\infty] := \Norm{\D u}_{\leb{\infty}(\W)}
\end{equation}
\cite{Aronsson:1965}
and as the (weighted) formal limit of the variational $p$-Laplacian
\begin{equation}
  \label{eq:introplap}
  \Delta_p u := \div{\qp{\norm{\D u}^{p-2} \D u}} = 0.
\end{equation}
The $p$-Laplacian is a divergence form problem and appropriate \emph{weak}
solutions to this problem are defined in terms of duality, or integration by parts.  In passing
to the limit ($p\to\infty$) the problem loses its divergence
structure. In the nondivergence setting we do not have access to the
same stability framework as in the variational case and a different
class of ``weak'' solution must be sought. The correct concept to use
is that of viscosity solutions
\cite[c.f.]{CrandallIshiiLions:1992,Katzourakis:2015}.
The main idea behind this solution concept is to pass derivatives to test functions through the maximum principle, that is, \emph{without} using duality.

The design of numerical schemes to approximate this solution concept
is limited, particularly in the finite element context, where the only
provably convergent scheme is given in \cite{JensenSmears:2013}
(although it is inapplicable to the problem at hand). In the finite
difference setting techniques have been developed
\cite{Oberman:2005,Oberman:2013} and applied to this problem and also
the associated eigenvalue problem \cite{Bozorgnia:2015}. In fact both
in the finite difference and finite element setting the methods of
convergence are based on the discrete monotonicity arguments of
\cite{BarlesSouganidis:1991} which is an extremely versatile
framework. Other methods exist for the problem, for example in
\cite{FengNeilan:2009}, the authors propose a biharmonic
regularisation which yields convergence in the case (\ref{eq:inflap})
admits a strong solution. In \cite{Pryer:2013} the author proposed an
$h$-adaptive finite element scheme based on a residual type error
indicator. The underlying scheme was based on the method derived in
\cite{LakkisPryer:2013} for fully nonlinear PDEs.

In this note we examine a different route. We will review and use the
known theory used in the derivation of the $\infty$-Laplacian
\cite[c.f.]{Aronsson:1986,Jensen:1993,Katzourakis:2015} where a
\emph{$p$-limiting} process is employed to derive
(\ref{eq:inflap}). We study how well Galerkin approximations of
(\ref{eq:introplap}) approximate the solutions of (\ref{eq:inflap})
and show that by forming an appropriate limit we are able to select
candidates for numerical approximation along a ``good'' sequence of
solutions. This is due to the equivalence of weak and viscosity
solutions to (\ref{eq:introplap})
\cite{JuutinenLindqvistManfredi:2001}. To be very clear about where
the novelty lies in this presentation, the techniques we use are not
new. We are summarising existing tools from two fields, one set from
PDE theory and the other from numerical analysis. While both sets of
results are relatively standard in their own field, to the authors'
knowledge, they have yet to be combined in this fashion.

We use this exposition to conduct some numerical experiments which
demonstrate the rate of convergence both in terms of $p$ approximation
\footnote{The terminology $p$ approximation we use here should not be
  confused with $p$-adaptivity which is local polynomial enrichment of
  the underlying discrete function space. } and
$h$ approximation. These results illustrate that for practical
purposes, as one would expect, the approximation of $p$-harmonic functions for large $p$
gives good resolution of $\infty$-harmonic functions. The numerical
approximation of $p$-harmonic functions is by now quite standard in
finite element literature, see for example \cite[\S
5.3]{Ciarlet:1978}. There has been a lot of activity in the area since
then however. In particular, the quasi-norm introduced in
\cite{BarrettLiu:1994} gave significant insight in the numerical
analysis of this problem and spawned much subsequent research for
which \cite{LiuYan:2001,Carstensen:2006,DieningKreuzer:2008} form an
inexhaustive list.

While it is not the focus of this work, we are interested in this
approach as it allows us to extend quite simply and reliably into the
vectorial case. When moving from scalar to vectorial calculus of
variations in $\leb{\infty}$ viscosity solutions are no longer
applicable. One notion of solution currently being investigated is
$\mathcal D$-solutions \cite{Katzourakis:2015a} which is based on
concepts of Young measures. The ultimate goal of this line of research
is the construction of reliable numerical schemes which allow for
various conjectures to be made as to the nature of solutions and even
what the correct solution concepts in the vectorial case are
\cite{KatzourakisPryer:2015}.

The rest of the paper is set out as follows: In \S \ref{sec:plap} we
formalise notation and begin exploring some of the properties of the
$p$-Laplacian. In particular, we recall that the notion of weak and
viscosity solutions to this problem coincide, allowing the passage to
the limit $p\to\infty$. In \S \ref{sec:fem} we describe a conforming
discretisation of the $p$-Laplacian and its properties. We show that
the method converges to the weak solution for fixed $p$. 
Numerical
experiments are given in \S \ref{sec:numerics} illustrating the
behaviour of numerical approximations to this problem.

\section{Approximation via the $p$-Laplacian}
\label{sec:plap}

In this section we describe how $\infty$-harmonic functions
can be approximated using $p$-harmonic functions. We give a brief
introduction to the $p$--Laplacian problem, beginning by introducing the Sobolev spaces
\cite{Ciarlet:1978,Evans:1998}
\begin{gather}
  \leb{p}(\W)
  =
  \ensemble{\phi}
           {\int_\W \norm{\phi}^p < \infty}   \text{ for } p\in[1,\infty) 
             \text{ and }
             \leb{\infty}(\W)
             =
             \ensemble{\phi}
                      {\esssup_\W \norm{\phi} < \infty},
                      \\
                      \sob{l}{p}(\W) 
                      = 
                      \ensemble{\phi\in\leb{p}(\W)}
                               {\D^{\vec\alpha}\phi\in\leb{p}(\W), \text{ for } \norm{\geovec\alpha}\leq l}
                               \text{ and }
    \sobh{l}(\W)
    := 
    \sob{l}{2}(\W),
\end{gather}
which are equipped with the following norms and semi-norms:
\begin{gather}
  \Norm{v}_{\leb{p}(\W)}^p
  :=
       {\int_\W \norm{v}^p} \text{ for } p \in [1,\infty)
         \text{ and }
  \Norm{v}_{\leb{\infty}(\W)}
  :=
  \esssup_\W |v|
       \\
       \Norm{v}_{l,p}^p
       := 
       \Norm{v}_{\sob{l}{p}(\W)}^p
       = 
       \sum_{\norm{\vec \alpha}\leq l}\Norm{\D^{\vec \alpha} v}_{\leb{p}(\W)}^p
       \\
       \norm{v}_{l,p}^p
       :=
       \norm{v}_{\sob{l}{p}(\W)}^p
       =
       \sum_{\norm{\vec \alpha} = l}\Norm{\D^{\vec \alpha} v}_{\leb{p}(\W)}^p
\end{gather}
where $\vec\alpha = \{ \alpha_1,\dots,\alpha_d\}$ is a
multi-index, $\norm{\vec\alpha} = \sum_{i=1}^d\alpha_i$ and
derivatives $\D^{\vec\alpha}$ are understood in a weak sense. We pay
particular attention to the case $l = 1$ and
\begin{gather}
  \sobg{1}{p}(\W) 
  :=
  \ensemble{\phi\in\sob{1}{p}(\W)}{\phi\vert_{\partial \W} = g},
\end{gather}
for a prescribed function $g\in\sob{1}{\infty}(\W)$.
Let $L = L\qp{\geovec x, u, \D u}$ be the \emph{Lagrangian}. We
will let
\begin{equation}
  \label{eq:action-functional}
  \dfunkmapsto[]
  {\cJ[\ \cdot \ ; p]}
  \phi
  {\sobg{1}{p}(\W)}
  {\cJ[\phi; p] := \int_\W L(\geovec x, \phi, \D \phi) \d \geovec x}
  {\reals}
\end{equation}
be known as the \emph{action functional}. For the $p$--Laplacian the
action functional is given as
\begin{equation}
  \cJ[u;p]
  :=
  \int_\W L(\geovec x, u, \D u)
  =
  \int_\W \norm{\D u}^p. \quad \footnote{Typically $L(\geovec x, u , \D u) = \tfrac{1}{p} \norm{\D u}^p$. Note here the rescaling of $L$ has no effect on the resultant Euler--Lagrange equations as to $L$ is independent of $u$.}
\end{equation}
We then look to find a
minimiser over the space $\sobg{1}{p}(\W)$, that is, to find
$u\in\sobg{1}{p}(\W)$ such that
\begin{equation}
  \cJ[u;p] = \min_{v\in\sobg{1}{p}(\W)} \cJ[v;p].
\end{equation}

If we assume temporarily that we have access to a smooth minimiser,
\ie $u\in\cont{2}(\W)$, then, given that the Lagrangian is of first order,
we have that the Euler--Lagrange equations are (in general) second
order.

The Euler--Lagrange equations for this problem are
\begin{equation}
  \label{eq:p-biharm}
  \cL[u; p] :=
  \div\qp{\norm{\D u}^{p-2}\D u} = 0.
\end{equation}
Note that, for $p=2$, the problem coincides with the Poisson
problem $\Delta u = 0$. In general, the $p$-Laplace problem is to find $u$ such that
\begin{equation}
  \label{eq:plap}
  \begin{split}
    \Delta_p u:= \div\qp{\norm{\D u}^{p-2} \D u} &= 0 \text{ in } \W
    \\
    u &= g \text{ on } \partial\W.
  \end{split}
\end{equation}

\begin{Defn}[weak solution]
  The problem (\ref{eq:plap}) is associated to a weak formulation, set
  \begin{gather}
    \bi{u}{v} = \int_\W \qp{\norm{\D u}^{p-2}\D u} \cdot \D v
  \end{gather}
  to be a semilinear form, then $u \in
  \sobg{1}{p}(\W)$ is a \emph{weak solution} of (\ref{eq:plap}) if it satisfies
  \begin{equation}
    \bi{u}{v} = 0 \Foreach v \in \sobz{1}{p}(\W).
  \end{equation}
\end{Defn}

\begin{Pro}[existence and uniqueness of weak solutions to (\ref{eq:plap})]
  \label{pro:existence-uniqueness}
  There exists a unique weak solution to (\ref{eq:plap}).
\end{Pro}
\begin{Proof}
  The proof is standard and can be found in \cite[Thm
  5.3.1]{Ciarlet:1978} for example. It is based on the strict
  convexity of $\cJ[\cdot; p]$, yielding uniqueness, together with
  appropriate growth conditions for existence.
\end{Proof}

\begin{Defn}[viscosity super and sub-solutions]
  \label{def:visc-soln}
  A function $u\in\cont{0}(\W)$ is a viscosity sub-solution of
  a general second order PDE 
  \begin{equation}
    \label{eq:fnl}
    F(u,\D u,\D^2 u) = 0
  \end{equation}
  at a point $\geovec x \in \W$ if for any $\phi\in\cont{2}(\W)$
  satisfying $u(\geovec x) = \phi(\geovec x)$, and touching $u$ from
  above, that is $u - \phi \leq 0$ in a neighbourhood of $\geovec x$,
  we have
  \begin{equation}
      F(\phi,\D \phi,\D^2 \phi) \leq 0.
  \end{equation}

  Similarly, a function $u\in\cont{0}(\W)$ is a viscosity
  super-solution of (\ref{eq:fnl}) at a point $\geovec x\in\W$ if for
  any $\phi\in\cont{2}(\W)$ satisfying $u(\geovec x) = \phi(\geovec
  x)$ and touches $u$ from below we have
  \begin{equation}
      F(\phi,\D \phi,\D^2 \phi) \geq 0.
  \end{equation}
\end{Defn}

\begin{Defn}[viscosity solution]
  The function $u$ is a viscosity solution of (\ref{eq:fnl}) in $\W$
  if it is a viscosity super and sub-solution at any $\geovec x
  \in\W$.
\end{Defn}

\begin{The}[weak solutions of the $p$-Laplacian are viscosity solutions]
  \label{the:weak-visc}
  Let $g\in\sob{1}{\infty}(\W)$ and suppose $p > d \geq 2$ is
  fixed. Then weak solutions of
  \begin{equation}
    \label{eq:p-lap}
    \begin{split}
      \Delta_p u &= 0 \text{ in } \W
      \\
      u &= g \text{ on } \partial\W
    \end{split}
  \end{equation}
  are viscosity solutions of
  \begin{equation}
    \label{eq:normalised}
    \begin{split}
      \Delta_\infty u + \frac{\norm{\D u}^2}{p-2} \Delta u &= 0 \text{ in } \W
      \\
      u &= g \text{ on } \partial\W.
    \end{split}
  \end{equation}
\end{The}
\begin{Proof}
  We begin by noting that by expanding the derivatives
  \begin{equation}
    \label{eq:norm}
    \begin{split}
      \Delta_p u &= \div\qp{\norm{\D u}^{p-2} \D u}
      \\
      &= 
      \norm{\D u}^{p-2} \Delta u
      + 
      \qp{p-2}\norm{\D u}^{p-4} \frob{\D u \otimes \D u}{\Hess u}
      \\
      &= 
      \norm{\D u}^{p-2} \Delta u
      + 
      \qp{p-2}\norm{\D u}^{p-4} \Delta_\infty u,
    \end{split}
  \end{equation}
  is a renormalisation of (\ref{eq:normalised}). The two formulations (\ref{eq:norm}) and (\ref{eq:normalised})
  of the $p$-Laplacian are equivalent in the viscosity sense, see for
  example \cite[\S 8 Lemma 3]{Katzourakis:2015}.

  It remains to show that weak solutions of (\ref{eq:p-lap}) are viscosity solutions of (\ref{eq:norm}). As $g\in\sob{1}{\infty}(\W)$ we have
  \begin{equation}
    \Norm{\D g}_{\leb{p}(\W)}^p
    =
    \int_\W \norm{\D g}^p
    =
    \cJ[g;p]
    < \infty.
  \end{equation}
  Since $u$ solves (\ref{eq:p-lap}) weakly, it minimises the functional $\cJ[\cdot;p]$ and hence the minimiser must be of finite energy. In view of the existence and uniqueness of the minimisation problem from Proposition \ref{pro:existence-uniqueness} and Morrey's inequality, we may infer $u\in\cont{0,\alpha}(\W)$ and hence $u\in\cont{0}(\W)$.

  Now assume by contradiction that $u$ is not a viscosity subsolution of
  \begin{equation}
    \norm{\D u}^{p-2} \Delta u
    + 
    \qp{p-2}\norm{\D u}^{p-4} \Delta_\infty u = 0
  \end{equation}
  then by Definition \ref{def:visc-soln} we can find an $\geovec x\in\W$, a $\psi\in\cont{2}(\W)$ and an $r>0$ such that $u - \psi <0$ on $B(\geovec x,r)$, $\qp{u - \psi}(\geovec x) = 0$ and
  \begin{equation}
    \Delta_p \psi
    =
    \norm{\D \psi}^{p-2} \Delta \psi
    + 
    \qp{p-2}\norm{\D \psi}^{p-4} \Delta_\infty \psi \leq C < 0 \text{ in } B(\geovec x, r),
  \end{equation}
  for some $C>0$. Since $u - \psi$ has a strict maximum at $\geovec x$ we may find an $\epsilon > 0$ such that
  \begin{equation}
    \w := \ensemble{\geovec x}{u(\geovec x) - \qp{\psi(\geovec x) - \epsilon} > 0}
    \subset
    B(\geovec x,r).
  \end{equation}
  Hence
  \begin{equation}
    \begin{split}
      C \int_\w \qp{u - \qp{\psi - \epsilon}}
      &\leq
      \int_\w - \Delta_p \psi \qp{u - \qp{\psi - \epsilon}}
      \\
      &=
      \int_\w \qp{\D \psi}^{p-2} \D \psi \cdot \D \qp{u - \psi},
    \end{split}
  \end{equation}
  as $u = \psi - \epsilon$ on $\partial \w$. Now by the convexity of the Lagrangian $L(\geovec x, u, \D u) = \norm{\D u}^p$ we have
  \begin{equation}
    \norm{\D v}^{p-2} \D v \cdot \D \qp{w - v} \leq \norm{\D w}^p - \norm{\D v}^p,
  \end{equation}
  hence
  \begin{equation}
    0 \leq C \int_\w \qp{u - \qp{\psi - \epsilon}}
    \leq
    \norm{\D u}^p - \norm{\D \qp{\psi - \epsilon}}^p
    =
    \cJ[u;p] - \cJ[\psi-\epsilon ;p]
  \end{equation}
  and we see
  \begin{equation}
    \cJ[u;p] \leq \cJ[\psi - \epsilon ;p].
  \end{equation}
  This means that $\w = \emptyset$ and we have a contradiction. The complete proof in full generality for convex minimisation problems can be found in \cite{Katzourakis:2015}. See also \cite{Katzourakis:2015b} where the author extends the arguments of \cite{JuutinenLindqvistManfredi:2001} to singular PDEs.
\end{Proof}

\begin{Rem}[viscosity solutions of the $p$-Laplacian are weak solutions]
  The converse to Theorem \ref{the:weak-visc} is also true, thus weak and viscosity solutions are an equivalent concept for the $p$-Laplacian and its evolutionary relative. This has been shown in \cite{JuutinenLindqvistManfredi:2001}.
\end{Rem}

\begin{The}[the limit as $p\to\infty$]
  \label{the:ptoinf}
  Let $u_p \in\sobg{1}{p}(\W)$ denote a sequence of weak/viscosity
  solutions to the $p$-Laplacian then there exists a subsequence such
  that as $p\to\infty$ that sequence converges to a candidate
  $\infty$-harmonic function $u_\infty\in\sob{1}{\infty}(\W)$, that is,
  \begin{equation}
    u_{p_j} \to u_\infty \text{ in } C^{0}.
  \end{equation}
\end{The}
\begin{Proof}
  We denote $u_p \in \sobg{1}{p}(\W)$ as the weak solution of (\ref{eq:plap}). In view of Proposition \ref{pro:existence-uniqueness} we know that $u_p$ minimises the energy functional
  \begin{equation}
    \cJ[u_p;p] = \int_\W \norm{\D u_p}^p.
  \end{equation}
  Hence in particular
  \begin{equation}
    \cJ[u_p;p] \leq \cJ[g;p],
  \end{equation}
  where $g$ is the associated boundary data to (\ref{eq:plap}). Using this fact
  \begin{equation}
    \Norm{\D u_p}_{\leb{p}(\W)}^p
    =
    \cJ[u_p;p]
    \leq
    \cJ[g;p]
    =
    \Norm{\D g}_{\leb{p}(\W)}^p,
  \end{equation}
  and we may infer that
  \begin{equation}
    \label{eq:limitpf1}
    \Norm{\D u_p}_{\leb{p}(\W)}
    \leq
    \Norm{\D g}_{\leb{p}(\W)}.
  \end{equation}
  Now fix a $k > d$ and take $p \geq k$, then using H\"olders inequality
  \begin{equation}
    \label{eq:limitpf3}
    \Norm{\D u_p}_{\leb{k}(\W)}^k
    =
    \int_\W \norm{\D u_p}^k
    \leq
    \qp{\int_\W 1^q}^{1/q} \qp{\int_\W \norm{\D u_p}^p}^{1/r},
  \end{equation}
  with $r = \tfrac{p}{k}$ and $q = \tfrac{r-1}{r}$ such that $\tfrac{1}{r} + \tfrac{1}{q} = 1$. Hence
  \begin{equation}
    \Norm{\D u_p}_{\leb{k}(\W)}^k
    \leq
    \norm{\W}^{\tfrac{r}{r-1}} \Norm{\D u_p}^{k}_{\leb{p}(\W)}
    =
    \norm{\W}^{1-\tfrac{k}{p}} \Norm{\D u_p}^{k}_{\leb{p}(\W)}
  \end{equation}
  and we see
  \begin{equation}
    \label{eq:limitpf2}
    \Norm{\D u_p}_{\leb{k}(\W)}
    \leq
    \norm{\W}^{\tfrac{1}{k}-\tfrac{1}{p}} \Norm{\D u_p}_{\leb{p}(\W)}.    
  \end{equation}
  Using the triangle inequality
  \begin{equation}
    \begin{split}
    \Norm{u_p}_{\leb{k}(\W)}
    &\leq
    \Norm{u_p - g}_{\leb{k}(\W)}
    +
    \Norm{g}_{\leb{k}(\W)}
    \\
    &\leq
    C_P \Norm{\D u_p - \D g}_{\leb{k}(\W)}
    +
    \Norm{g}_{\leb{k}(\W)},
    \end{split}
  \end{equation}
  in view of the Poincar\'e inequality. Using the triangle inequality again we have
  \begin{equation}
    \begin{split}
    \Norm{u_p}_{\leb{k}(\W)}
    &\leq
    C_P \Norm{\D u_p}_{\leb{k}(\W)}
    +
    \Norm{g}_{\sob{1}{k}(\W)}
    \\
    &\leq
    C_P \norm{\W}^{\tfrac{1}{k}-\tfrac{1}{p}} \Norm{\D u_p}_{\leb{k}(\W)}
    +
    \Norm{g}_{\sob{1}{k}(\W)},
    \end{split}
  \end{equation}
  by (\ref{eq:limitpf2}). Hence using (\ref{eq:limitpf1})
  \begin{equation}
    \label{eq:limitpf4}
    \begin{split}
      \Norm{u_p}_{\sob{1}{k}(\W)}
      &\leq
      C \Norm{g}_{\sob{1}{k}(\W)}.
    \end{split}
  \end{equation}
  This means that for any $k > d$ we have uniformly that
  \begin{equation}
    \sup_{p > k} \Norm{u_p}_{\sob{1}{k}(\W)} \leq C.
  \end{equation}
  Hence, in view of weak compactness, we may extract a subsequence $\{
  u_{p_j} \}_{j=1}^\infty \subset \{ u_{p} \}_{p=1}^\infty$ and a
  function $u_\infty\in\sob{1}{k}(\W)$ such that for any $k > n$
  \begin{equation}
    u_{p_j} \rightharpoonup u_\infty \text{ weakly in } \sob{1}{k}(\W)
  \end{equation}
  and
  \begin{equation}
    \begin{split}
      \Norm{u_\infty}_{\sob{1}{k}(\W)} &\leq \liminf_{j\to\infty} \Norm{u_{p_j}}_{\sob{1}{k}(\W)}
      \\
      &\leq
      \liminf_{j\to\infty} 
      C\Norm{g}_{\sob{1}{k}(\W)}.
    \end{split}
  \end{equation}
  Taking the limit $k\to\infty$ we have
  \begin{equation}
    \begin{split}
      \Norm{u_\infty}_{\sob{1}{\infty}(\W)} &\leq C \Norm{g}_{\sob{1}{\infty}(\W)}
    \end{split}
  \end{equation}
  and thus
  $u_\infty \in \sob{1}{\infty}(\W)$. The result follows from Morrey's
  inequality, concluding the proof. 
\end{Proof}

\begin{Rem}[An alternative to the $p$-Dirichlet functional]
  We note that an alternative sequence of solutions is given in \cite{EvansSmart:2011}, where rather than studying the limit of the $p$-Dirichlet functional, the authors propose 
  \begin{equation}
    \widetilde{\cJ}[u;p]
    =
    \int_\W \exp\qp{{p}\norm{\D u}^2}.
  \end{equation}
  This functional may have some merit over the $p$-Dirichlet functional since the Euler--Lagrange equations
  \begin{equation}
    \begin{split}
      0 &= \div\qp{\exp\qp{p\norm{\D u}^2}\D u}
      \\
      &= \exp\qp{p\norm{\D u}^2} \Delta u + p \exp\qp{p\norm{\D u}^2} \Delta_\infty u
    \end{split}
  \end{equation}
  yield a clearer relation between $\Delta u$ and $\Delta_\infty
  u$. We will not explore this issue further in this work.
\end{Rem}

\begin{The}[existence and uniqueness of viscosity solutions to the
  $\infty$-Laplacian \cite{Jensen:1993}]
  The ``candidate'' $\infty$-harmonic function $u_\infty$ is the
  unique viscosity solution to the $\infty$-Laplacian
  (\ref{eq:inflap}).
\end{The}
\begin{Proof}
  The proof is detailed in \cite{Jensen:1993}. Roughly, existence of
  $u_\infty$ has been shown in Theorem \ref{the:ptoinf}. For
  uniqueness one must prove and make use of the maximum principle for (\ref{eq:inflap}). Note also the result of \cite{ArmstrongSmart:2010} where the authors use difference equations to prove the same result in a simpler fashion.
\end{Proof}

\section{Discretisation of the $p$-Laplacian}
\label{sec:fem}

In this section we describe a conforming finite element discretisation of the $p$-Laplacian.
Let $\T{}$ be a conforming triangulation of $\W$,
namely, $\T{}$ is a finite family of sets such that
\begin{enumerate}
\item $K\in\T{}$ implies $K$ is an open simplex (segment for $d=1$,
  triangle for $d=2$, tetrahedron for $d=3$),
\item for any $K,J\in\T{}$ we have that $\closure K\meet\closure J$ is
  a full lower-dimensional simplex (i.e., it is either $\emptyset$, a vertex, an
  edge, a face, or the whole of $\closure K$ and $\closure J$) of both
  $\closure K$ and $\closure J$ and
\item $\union{K\in\T{}}\closure K=\closure\W$.
\end{enumerate}
The shape regularity constant of $\T{}$ is defined as the number
\begin{equation}
  \label{eqn:def:shape-regularity}
  \mu(\T{}) := \inf_{K\in\T{}} \frac{\rho_K}{h_K},
\end{equation}
where $\rho_K$ is the radius of the largest ball contained inside
$K$ and $h_K$ is the diameter of $K$. An indexed family of
triangulations $\setof{\T n}_n$ is called \emph{shape regular} if 
\begin{equation}
  \label{eqn:def:family-shape-regularity}
  \mu:=\inf_n\mu(\T n)>0.
\end{equation}
Further, we define $\funk h\W\reals$ to be the {piecewise
  constant} \emph{meshsize function} of $\T{}$ given by
\begin{equation}
h\equiv  h(\vec{x}):=\max_{\closure K\ni \vec{x}}h_K.
\end{equation}
A mesh is called quasiuniform when there exists a positive constant
$C$ such that $\max_{x\in\Omega} h \le C \min_{x\in\Omega} h$. In what
follows we shall assume that all triangulations are shape-regular and
quasiuniform although the results may be extendable even in the
non-quasiuniform case using techniques developed in
\cite{DieningKreuzer:2008}.

We let $\E{}$ be the skeleton (set of common interfaces) of the
triangulation $\T{}$ and say $e\in\E$ if $e$ is on the interior of
$\W$ and $e\in\partial\W$ if $e$ lies on the boundary $\partial\W$
{and set $h_e$ to be the diameter of $e$.}

Further, we define the broken gradient $\D_h$, Laplacian
$\Delta_h$ and Hessian $\Hess_h$ to be defined element-wise by
$\D_h w|_K = \D w$, $\Delta_h w|_K = \Delta w$, $\Hess_h w|_K
= \Hess w$ for all $K\in \T{}$, respectively, for respectively smooth
functions on the interior of $K$,

We let $\poly k(\T{})$ denote the space of piecewise polynomials of
degree $k$ over the triangulation $\T{}$,\ie
\begin{equation}
  \poly k (\T{}) = \{ \phi \text{ such that } \phi|_K \in \poly k (K) \}
\end{equation}
 and introduce the \emph{finite element space}
\begin{gather}
  \label{eqn:def:finite-element-space}
  \fes := \poly k(\T{}) \cap \cont{0}(\W)
\end{gather}
to be the usual space of continuous piecewise polynomial
functions of degree $k$.


\begin{Defn}[finite element sequence]
  A finite element sequence $\{ V, \fes \}$ is a sequence of
  discrete objects indexed by the mesh parameter, $h$, and
  individually represented on a particular finite element space
  $\fes$, which itself has a discretisation parameter $h$, that is
  $\fes = \fes(h)$.
\end{Defn}

\begin{Defn}[$\leb{2}(\W)$ projection operator]
  \label{eq:l2proj}
  The $\leb{2}(\W)$ projection operator, $P_k : \leb{2}(\W) \to \fes$ is
  defined for $v\in\leb{2}(\W)$ such that
  \begin{equation}
    \int_\W P_k v \Phi = \int_\W v \Phi \Foreach \phi\in\fes.
  \end{equation}
  It is well known that this operator satisfies the following
  approximation properties for $v \in \sob{1}{p}(\W)$ 
  \begin{gather}
    \lim_{h\to 0}\Norm{v - P_k v}_{\leb{p}(\W)} = 0
    \\
    \lim_{h\to 0}\Norm{\D v - \D \qp{ P_k v}}_{\leb{p}(\W)} = 0.
  \end{gather}
\end{Defn}

\subsection{Galerkin discretisation}
We consider the Galerkin discretisation of (\ref{eq:plap}), to find $U \in \fes$ with $U\vert_{\partial \W} = P_k g$ such that
\begin{equation}
  \label{eq:plapdis}
  \bi{U}{\Phi} = 0 \Foreach \Phi \in \fes.
\end{equation}

\begin{Pro}[existence and uniqueness of solution to (\ref{eq:plapdis})]
  There exists a unique solution of (\ref{eq:plapdis}).
\end{Pro}
\begin{Proof}
  The proof is standard and, in fact, equivalent to that of the smooth case, as in \cite[Thm 5.3.1]{Ciarlet:1978}.
\end{Proof}

\begin{The}[convergence of the discrete scheme to weak solutions]
  \label{the:convergence-of-scheme}
  Let $\{ U_p, \fes\}$ be the finite element sequence generated by
  solving (\ref{eq:plapdis}) and $u_p$, the weak solution of (\ref{eq:p-lap}), then for fixed $p$ we have that
  \begin{equation}
    U_p \to u_p  \text{ in } \cont{0}(\W).
  \end{equation}
\end{The}
\begin{Proof}
  We begin by noting the discrete weak formulation (\ref{eq:plapdis})
  is equivalent to the minimisation problem: Find $U\in\fes$ such that
  \begin{equation}
    \label{eq:dismin}
    \cJ[U; p] = \min_{V\in\fes} \cJ[V; p].
  \end{equation}
  Using this, we immediately have
  \begin{equation}
    \Norm{\D U}_{\leb{p}(\W)}^p \leq \cJ[U; p] \leq \cJ[P g; p] \leq \Norm{\D \qp{P g}}_{\leb{p}(\W)}^p.
  \end{equation}
  In view of the stability of the $\leb{2}$ projection in
  $\sob{1}{p}(\W)$ \cite{CrouzeixThomee:1987} we have
  \begin{equation}
    \Norm{\D U}_{\leb{p}(\W)} \leq C,
  \end{equation}
  uniformly in $h$.
  Hence by weak compactness there exists a (weak) limit to the finite
  element sequence, which we will call $u^*$. Due to the weak
  semicontinuity of $\cJ[\cdot; p]$ we have 
  \begin{equation}
    \cJ[u^*; p] \leq \cJ[U; p].
  \end{equation}
  In addition, in view of the approximation properties of $P_k$ given in
  \ref{eq:l2proj} we have for any $v\in\cont{\infty}$ that
  \begin{equation}
    \cJ[v; p] = \liminf_{h\to 0} \cJ[P_k v; p].
  \end{equation}
  Using the fact that $U$ is a discrete minimiser of (\ref{eq:dismin})
  we have
  \begin{equation}
    \cJ[u^*; p] \leq \cJ[U; p] \leq \cJ[P_k v; p],
  \end{equation}
  whence sending $h\to 0$ we see 
  \begin{equation}
    \cJ[u^*;p] \leq \cJ[v;p].
  \end{equation}
  Now, as $v$ was generic we may use density arguments and that $u_p$ was the unique minimiser to conclude $u^* = u_p$, concluding the proof.
\end{Proof}

\begin{Rem}[convergence of the discrete scheme to viscosity solutions]
  In view of Theorem \ref{the:weak-visc} the discrete scheme converges
  to viscosity solutions of the $p$-Laplacian. 
\end{Rem}

\begin{Lem}[convergence in the limit $p\to\infty$]
  \label{eq:lem-dis-p-to-inf}
  Let $U_p$ solve the discrete problem (\ref{eq:plapdis}), then for fixed $h$ along
  a subsequence we have $U_p \to U_\infty$.
\end{Lem}
\begin{Proof}
  The proof follows similarly to Theorem \ref{the:ptoinf}. Since $U_p$ is the Galerkin solution to (\ref{eq:plapdis}), it minimises $\cJ[\cdot,p]$  over $\fes$. Hence we know
  \begin{equation}
    \Norm{\D U_p}_{\leb{p}(\W)}
    \leq
    \Norm{\D \qp{P_k g}}_{\leb{p}(\W)}
    \leq
    C \Norm{\D g}_{\leb{p}(\W)},
  \end{equation}
  in view of the stability of $P_k$ in $\sob{1}{p}(\W)$. In addition, analogously to (\ref{eq:limitpf3})--(\ref{eq:limitpf4}) we may find a constant such that
  \begin{equation}
    \Norm{U_p}_{\sob{1}{k}} \leq C,
  \end{equation}
  allowing the extraction of a subsequence $\{ U_{p_j}\}_{j=1}^\infty$ and a limit $U_\infty$ such that for $k > d$
  \begin{equation}
    U_{p_j} \rightharpoonup U_\infty \text{ weakly in }\sob{1}{k}(\W).
  \end{equation}
  The rest of the proof parallels that of Theorem \ref{the:ptoinf}.
\end{Proof}

\begin{Rem}[Summarising the results thus far]
  Up to this point we have shown the green (solid) lines on the following diagram
  \begin{figure}[h!]
  \begin{tikzpicture}
    \SetGraphUnit{5} 

\GraphInit[vstyle=Empty] 
\tikzset{VertexStyle/.append style = {shape=rectangle,inner sep=0pt}} 

\Vertex[L=$U_\infty$]{1} 
\EA[unit=4,L=$u_\infty$](1){2} 
\NO[unit=3,L=$U_p$](1){4} 
\NO[unit=3,L=$u_p$](2){3}

\begin{scope}[every node/.style={midway},>=latex']  
  \draw[->]        (4)--(3) node [above] {Thm \ref{the:convergence-of-scheme}};
  \draw[->]        (4)--(3) node [below] {$h \to 0$};
  \draw[->,green]  (4)--(3);
  \draw[->]        (3)--(2) node [right] {$p \to \infty$};
  \draw[->]        (3)--(2) node [left]  {Thm \ref{the:ptoinf}};
  \draw[->,green]  (3)--(2);
  \draw[->]        (4)--(1) node [left]  {$p \to \infty$};
  \draw[->]        (4)--(1) node [right]  {Lem \ref{eq:lem-dis-p-to-inf}};
  \draw[->,green]  (4)--(1);
  \draw[dashed,red](4)--(2); 
\end{scope}
  \end{tikzpicture}
  \end{figure}
  hold. We would like to select a route for which we can pass the limits together, that is, we want to select an appropriate route for which the red (dashed) line is true.
\end{Rem}

\begin{The}[convergence]
  \label{the:convergence}  
  Let $U_p$ be the Galerkin solution of (\ref{eq:plapdis}) and $u_\infty$ the unique viscosity solution of (\ref{eq:inflap}) then along a subsequence
  \begin{equation}
    U_{p_j} \to u_\infty \text{ in } \cont{0} \text{ as } p \to \infty \text{ and } h \to 0
  \end{equation}
\end{The}
\begin{Proof}
  The proof is a consequence of Theorems \ref{the:ptoinf} and \ref{the:convergence-of-scheme} noting that along the same subsequence used in Theorem \ref{the:ptoinf} we have that
  \begin{equation}
    \Norm{U_{p_j} - u_\infty}_{\cont{0}(\W)}
    \leq
    \Norm{U_{p_j} - u_{p_j}}_{\cont{0}(\W)}
    +
    \Norm{u_{p_j} - u_\infty}_{\cont{0}(\W)}
  \end{equation}
  and hence $\Norm{U_{p_j} - u_\infty}_{\cont{0}(\W)} \to 0$ as $p\to \infty$ and $h\to 0$.
\end{Proof}

\begin{Rem}[consequences of Theorem \ref{the:convergence}]
  An immediate consequence of Theorem \ref{the:convergence} and the
  previous arguments are that for $H : \W \times \reals \times
  \reals^d$ with appropriate conditions (convexity for example) ,
  finite element approximations to the $p$-functional
  \begin{equation}
    \cJ[u;p] = \Norm{H(\cdot, u, \D u)}_{\leb{p}(\W)}
  \end{equation}
  can be used as approximations to 
  \begin{equation}
    \cJ[u;\infty] = \Norm{H(\cdot, u, \D u)}_{\leb{\infty}(\W)}.
  \end{equation}
\end{Rem}


\begin{Rem}[discontinuous Galerkin approximations]
  All the above results can be extended into the discontinuous Galerkin framework. This is based on the discrete action functional
  \begin{equation}
    \label{eq:action-dg}
    \cJ_h[U; p] := \int_\W \norm{G(U)}^p + \int_\E h_e^{1-p} \norm{\jump{U}}^p,
  \end{equation}
  where
  \begin{equation}
    \int_\W G(U) \phi = \int_\W \D_h U \phi - \int_\E \jump{U} \avg{\phi} \Foreach \phi \in \poly{k}(\T{}),
  \end{equation}
  where $\jump{U} = U|_{K^+} - U|_{K^-}$ denotes the \emph{jump} over an edge $e$ shared by neighbouring elements $K^+$ and ${K^-}$ and $\avg{\phi} = \tfrac{1}{2} \qp{\phi|_{K^+} + \phi|_{K^-}}$, the \emph{average} of a quantity over an edge.
  Using the results of \cite{BurmanErn:2008} discrete minimisers to (\ref{eq:action-dg}) satisfy the equivalent weak convergence results to the conforming finite elements.
\end{Rem}

\section{Numerical experiments}
\label{sec:numerics}

In this section we summarise numerical experiments validating the analysis done in previous sections.

\begin{Rem}[practical computation of (\ref{eq:plapdis}) for large $p$]
  The computation of $p$-harmonic functions is an extremely challenging problem in its own right. The class of nonlinearity in the problem results in the algebraic system, which ultimately yields the finite element solution, being extremely badly conditioned. One method to tackle this class of problem is the use preconditioners based on descent algorithms \cite{HuangLiLiu:2007}. For extremely large $p$, say $p \geq 1000$ this may be required, however for our purposes we restrict our attention to $p \sim 100$. This yields sufficient accuracy for the results we want to illustrate.

  Even tackling the case $p\sim 100$ is computationally tough. Our numerical approximation is based on a Newton solver. As is well known, Newton solvers require a sufficiently close initial guess to converge. For large $p$ a reasonable initial guess is given by numerically approximating the $q$-Laplacian for $q < p$ sufficiently close to $p$. This leads to an iterative process in the generation of the initial guess.
\end{Rem}

\subsection{Test 1 : Approximation of the Aronsson viscosity solution}

\begin{figure}[!ht]
  \caption[Numerical Results for Problem \eqref{eqn:Problem:1} with
  $\poly1$ elements]
  {\label{fig:aronson}
    Finite element approximations to the $\infty$-harmonic Aronsson function \eqref{eq:aronson} using $p$-harmonic functions for various $p$. Notice as $p$ increases the approximation better catches the singularity on the coordinate axis.
  }
  \begin{center}
    \subfigure[{\label{fig:a1}
        The finite element approximation to the $5$-Laplacian.
    }]{
      \includegraphics[scale=\figscale,width=0.4\figwidth]{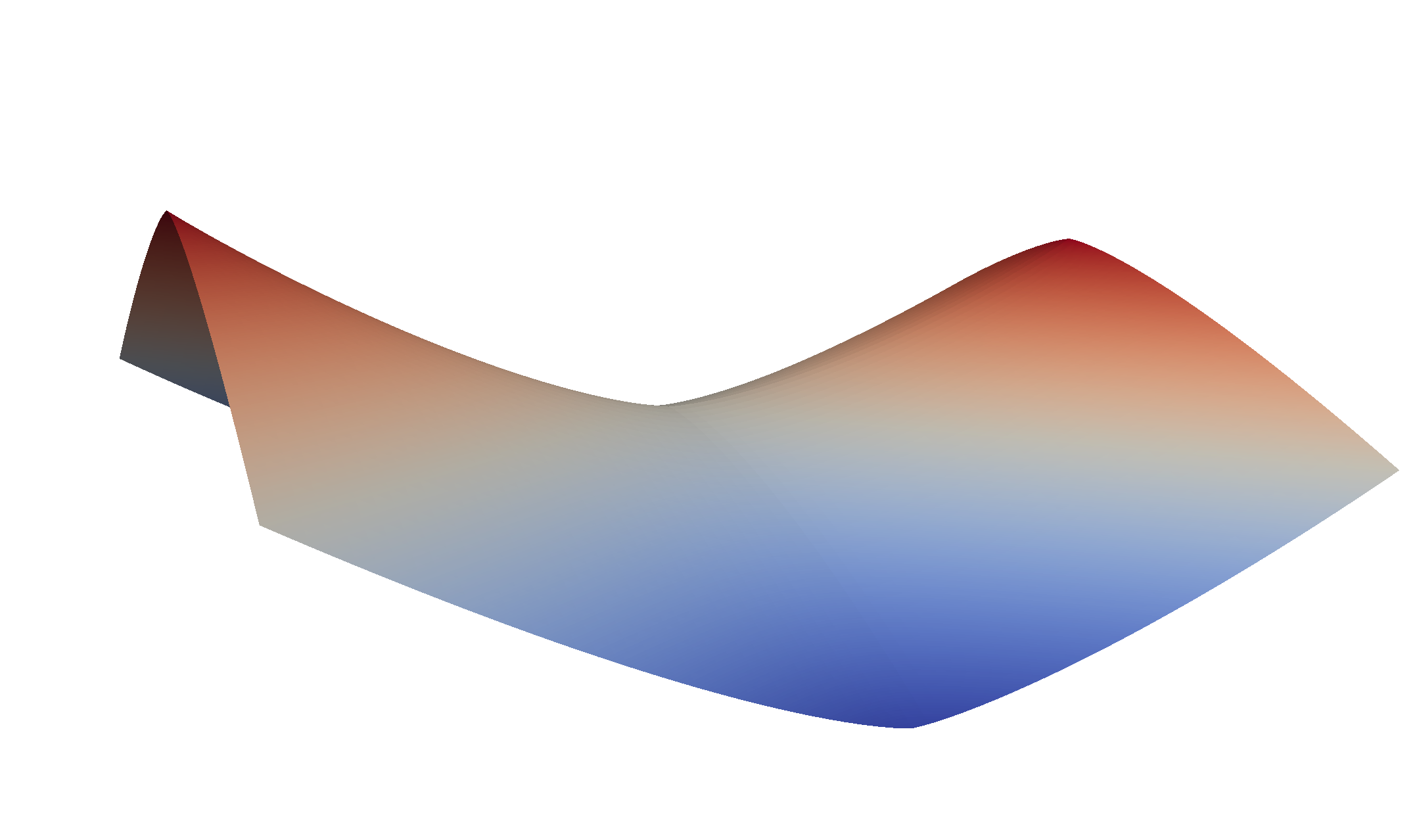}
    }
    \hfill
    \subfigure[{\label{fig:a2}
        The finite element approximation to the $15$-Laplacian.
    }]{
      \includegraphics[scale=\figscale,width=0.4\figwidth]{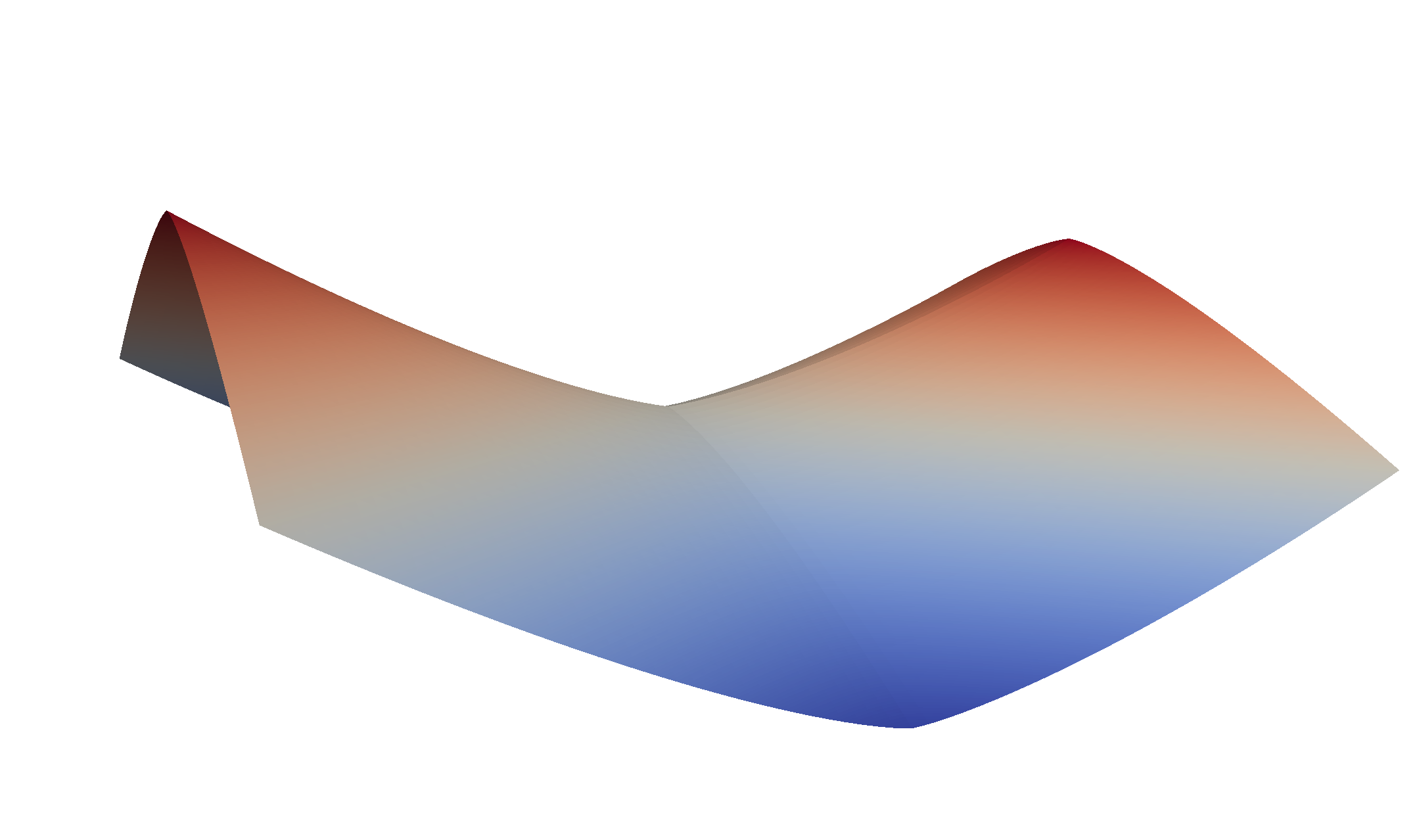}
    }
    \hfill
    \subfigure[{\label{fig:a2}
        The finite element approximation to the $50$-Laplacian.
    }]{
      \includegraphics[scale=\figscale,width=0.4\figwidth]{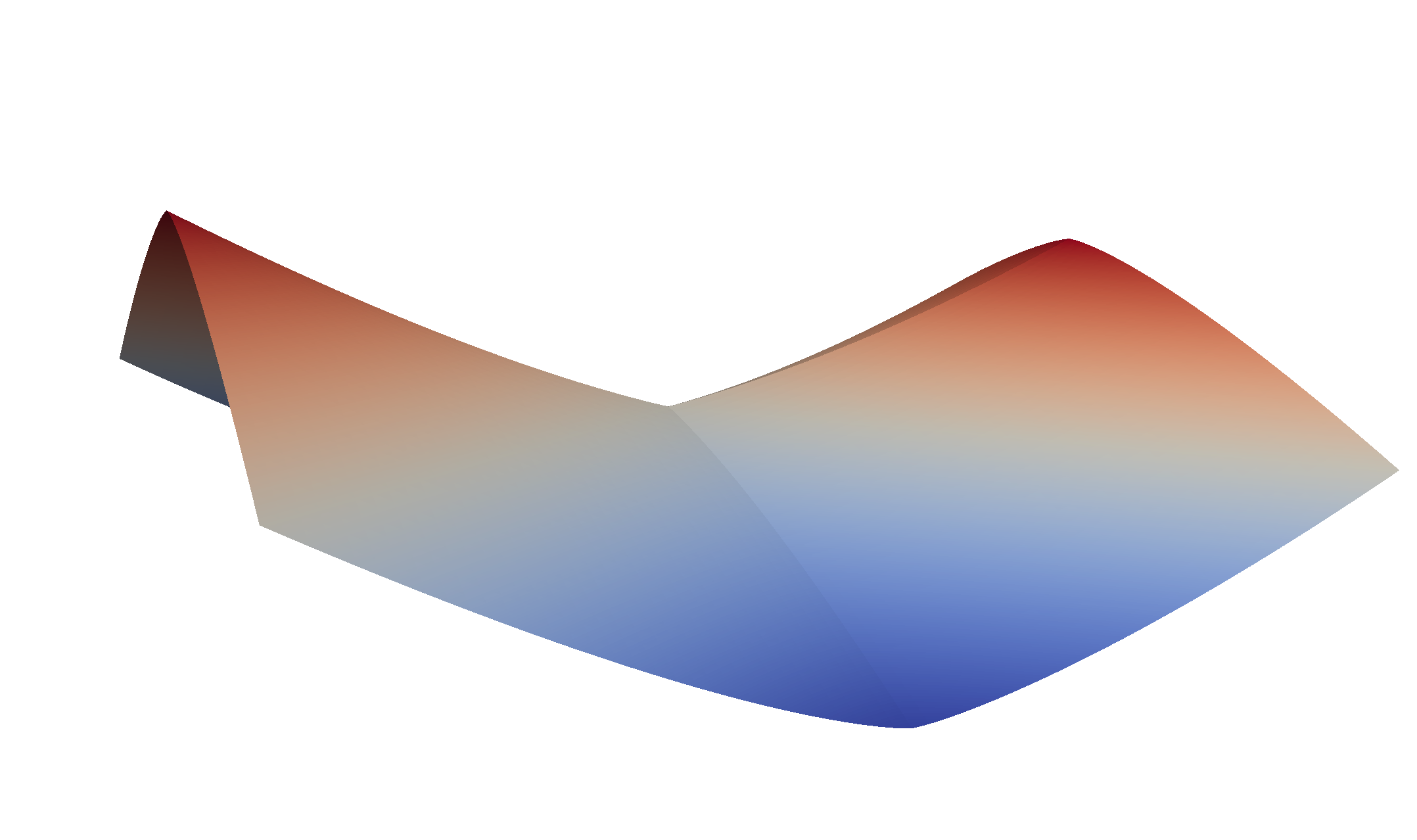}
    }
    \hfill
    \subfigure[{\label{fig:a2}
        The finite element approximation to the $100$-Laplacian.
    }]{
      \includegraphics[scale=\figscale,width=0.4\figwidth]{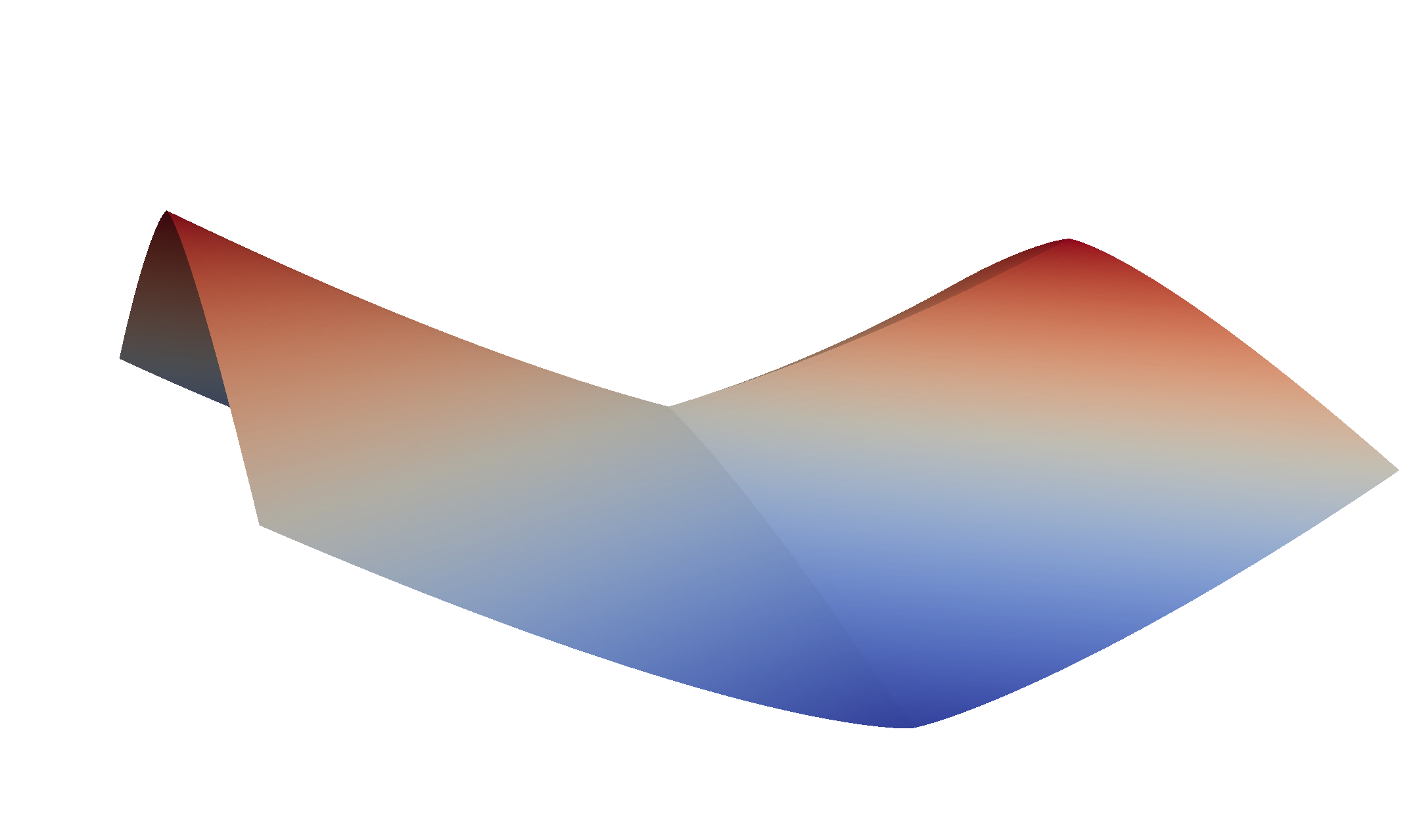}
    }
    \end{center}
  \end{figure}

\begin{figure}[ht]
    \caption{
    \label{fig:fixh}
    The error of the finite element approximation to the $p$-Laplacian compared to the viscosity solution to the $\infty$-Laplacian for various $p$. The colours represent different mesh refinement levels. The darker the colour, the more refined the mesh. In this experiment the meshsize ranges from $h\sim 0.7$ to $h \sim 0.005$ Notice as the mesh is refined, best approximation is achieved for higher and higher $p$.
  }
\centering
\begin{center}
  \includegraphics[scale=1.,width=.75\linewidth]{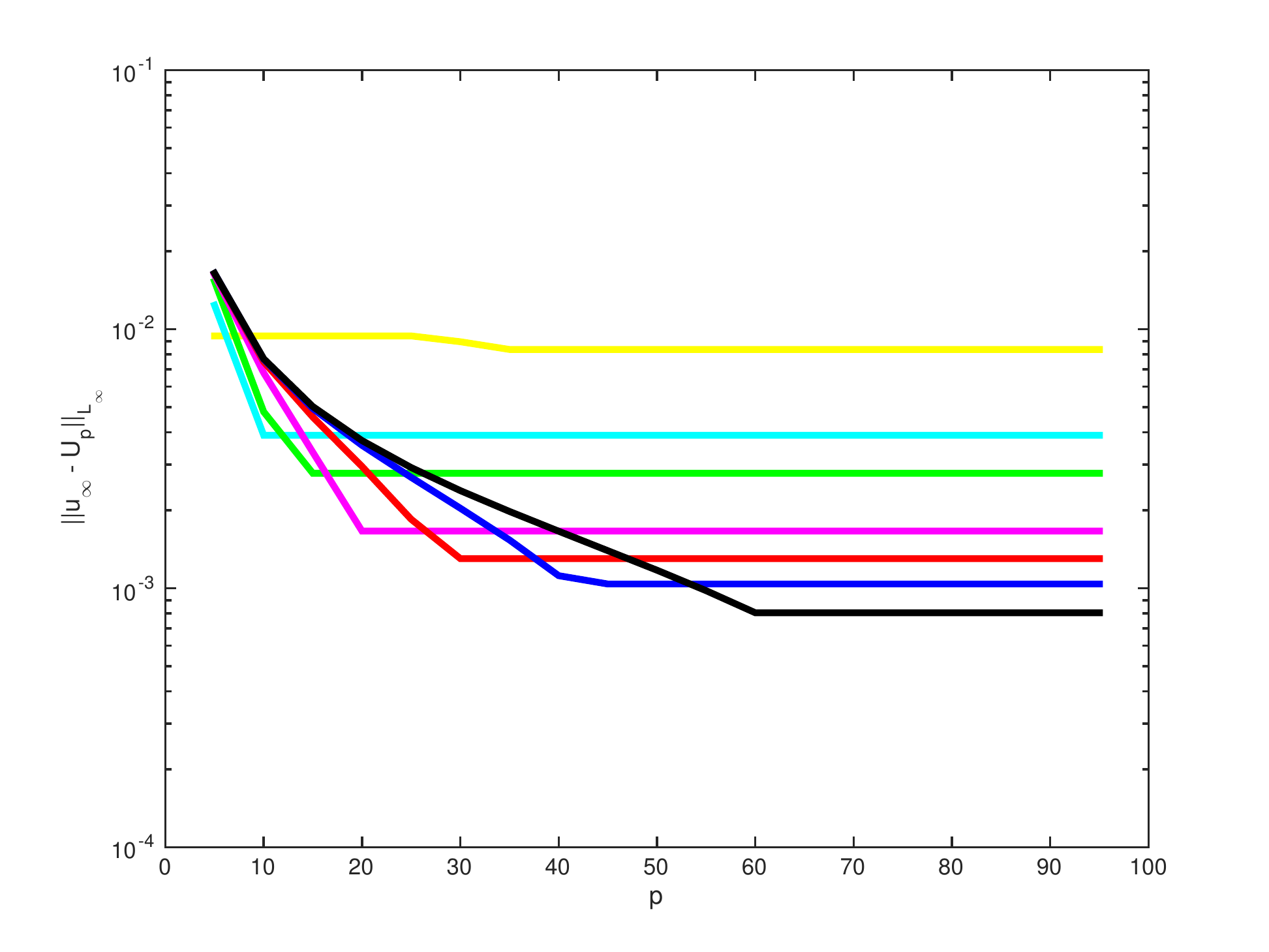}
\end{center}
\end{figure}

We begin by approximating the viscosity solution derived by Aronsson using separation of variables \cite{Aronsson:1986}. The function
\begin{equation}
  \label{eq:aronson}
  u(x,y) = \frac{3}{8}\qp{\norm{x}^{4/3} - \norm{y}^{4/3}} \in \cont{1,1/3}(\W)
\end{equation}
is a viscosity solution of the $\infty$-Laplacian. Notice that this is a weighted version of the Aronsson solution. We have chosen this as $\norm{\D u} \leq 1$ on the domains we consider to try to overcome the severe restrictions in computing $p$-harmonic functions with large $p$. In this test we take $\W = [-1.0001, 0.9999]^2$ and triangulate with a criss-cross mesh. This is so the singularity will not be aligned with the mesh. We approximate the solution of the $p$-Laplacian with boundary data given by (\ref{eq:aronson}) for a variety of increasing $p$. Examples of solutions are given in Figure \ref{fig:aronson}. In Figure \ref{fig:fixh} we plot the error against $p$ for a various levels of mesh refinement. In Table \ref{Table1} we demonstrate the convergence of the finite element approximations as $h\to 0$. 

\begin{table}[!ht]
  \caption{
    \label{Table1}
    In this Table we show the convergence of the finite element approximation $U_p$ to $u_\infty$, a viscosity solution of (\ref{eq:inflap}), as the meshsize is decreased. We study the $\leb{\infty}$ error of the approximation, the associated convergence rate and give $p^*$, the smallest such $p$ for which $\inf_p \Norm{u_\infty - U_p}_{\leb{\infty}(\W)}$ is attained. Notice that as the mesh is refined, the critical value increases.
   }
  \begin{center}
    \begin{tabular}{ | c | c | c | c | }
      \hline
      $\dim{\fes}$ & $\inf_p \Norm{u_\infty - U_p}_{\leb{\infty}(\W)}$ & EOC & $p^*$
      \\
      \hline
      $25$ & $0.0162$    & $0.00$ & $5$  \\
      \hline 
      $81$ & $0.00836$    & $0.95$ & $5$  \\
      \hline
      $289$ & $0.00390$    & $1.10$ & $10$  \\
      \hline
      $1089$ & $0.00278$    & $0.49$ & $15$  \\
      \hline
      $4225$ & $0.00166$    & $0.74$ & $20$  \\
      \hline
      $16641$ & $0.00130$    & $0.35$ & $30$  \\
      \hline
      $66049$ & $0.00104$    & $0.33$ & $45$  \\
      \hline
      $263169$ & $0.000805$    & $0.37$ & $60$  \\
      \hline
    \end{tabular}
  \end{center}
\end{table} 

\subsection{Test 2 : Approximation of a smooth solution}

\begin{figure}[ht]
    \caption{
    \label{fig:fixh2}
    The error of the finite element approximation to the $p$-Laplacian compared to the viscosity solution to the $\infty$-Laplacian for various $p$. The colours represent different mesh refinement levels. The darker the colour, the more refined the mesh. In this experiment the meshsize ranges from $h\sim 0.7$ to $h \sim 0.005$ Notice as the mesh is refined, best approximation is achieved for higher and higher $p$.
  }
\centering
\begin{center}
  \includegraphics[scale=1.,width=.75\linewidth]{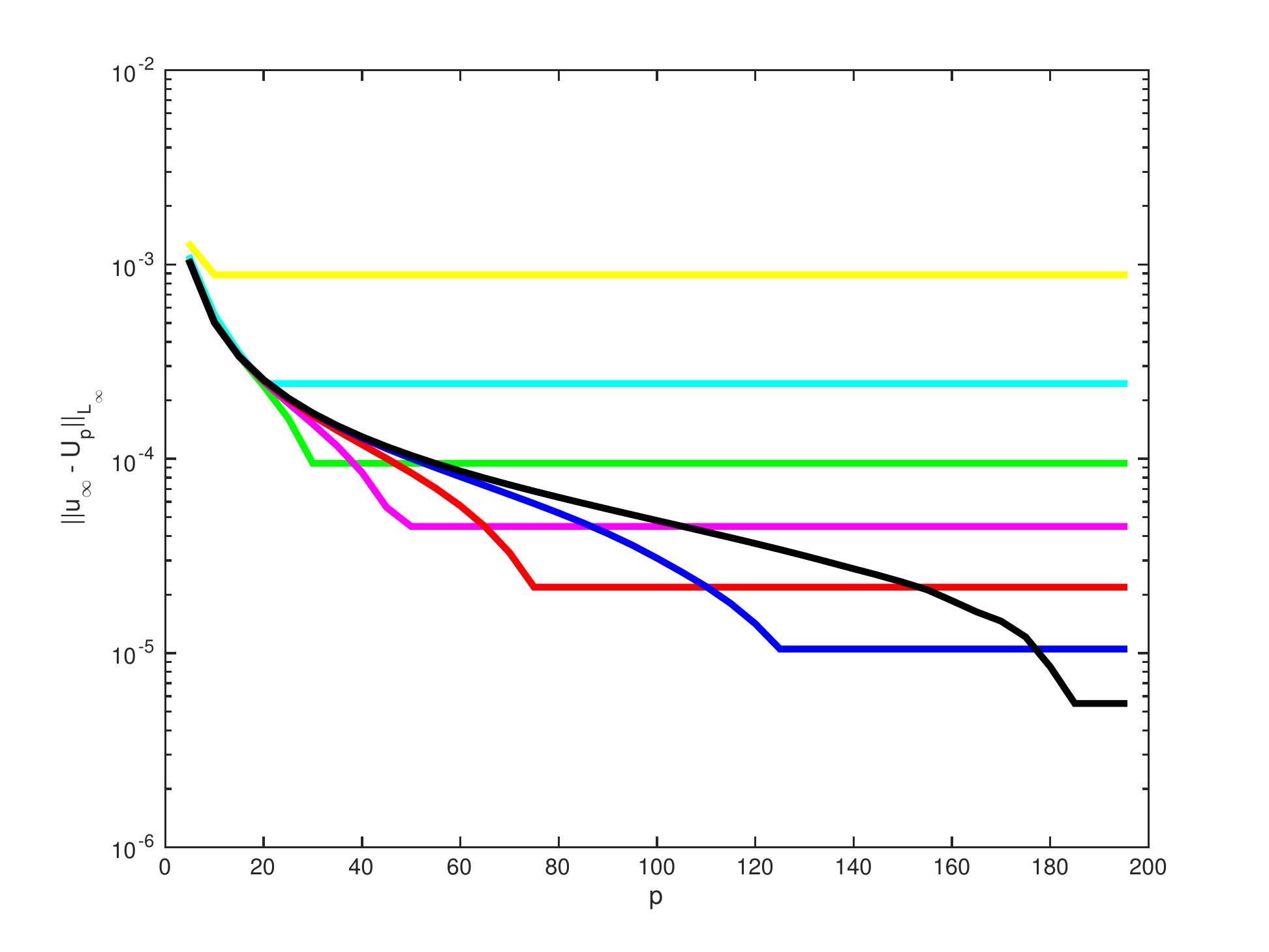}
\end{center}
\end{figure}

To test the approximation of a known smooth solution of (\ref{eq:inflap}) we look at (\ref{eq:aronson}) away from the coordinate axis. In this test we take $\W = [0.5, 1.5]^2$ and triangulate with a criss-cross mesh. As in Test 1, we approximate the solution of the $p$-Laplacian with boundary data given by (\ref{eq:aronson}) for a variety of increasing $p$. In Figure \ref{fig:fixh2} we plot the error against $p$ for a various levels of mesh refinement. In Table \ref{Table2} we demonstrate the convergence of the finite element approximations as $h\to 0$. 

\begin{table}[!ht]
  \caption{
    \label{Table2}
    In this Table we show the convergence of the finite element approximation $U_p$ to $u_\infty$, a smooth solution of (\ref{eq:inflap}), as the meshsize is decreased. We study the $\leb{\infty}$ error of the approximation, the associated convergence rate and give $p^*$, the smallest such $p$ for which $\inf_p \Norm{u_\infty - U_p}_{\leb{\infty}(\W)}$ is attained. Notice that as the mesh is refined, the critical value increases much quicker than in the nonsmooth case of Table \ref{Table1}.
  }
  \begin{center}
    \begin{tabular}{ | c | c | c | c | }
      \hline
      $\dim{\fes}$ & $\inf_p \Norm{u_\infty - U_p}$ & EOC & $p^*$
      \\
      \hline
      $25$ & $0.00301$    & $0.00$ & $5$  \\
      \hline 
      $81$ & $0.000883$    & $1.77$ & $10$  \\
      \hline
      $289$ & $0.000244$    & $1.86$ & $20$  \\
      \hline
      $1089$ & $0.0000946$    & $1.37$ & $30$  \\
      \hline
      $4225$ & $0.0000448$    & $1.08$ & $50$  \\
      \hline
      $16641$ & $0.0000218$    & $1.04$ & $75$  \\
      \hline
      $66049$ & $0.0000105$    & $1.06$ & $125$  \\
      \hline
      $263169$ & $0.0000052$    & $1.01$ & $185$  \\
      \hline
    \end{tabular}
  \end{center}
\end{table}

\end{document}